\documentclass[reqno,12pt]{amsart}

\usepackage{amssymb}
\usepackage{amsthm}

\newtheorem{theorem}{Theorem}[section]
\newtheorem{lemma}[theorem]{Lemma}

\theoremstyle{definition}

\numberwithin{equation}{section}

\begin{document}

\title[Asymptotic values of entire functions of infinite order]{Asymptotic values of entire functions of infinite order}

\author[Aimo Hinkkanen and Joseph Miles]{Aimo Hinkkanen$^{1}$, Joseph Miles$^{2}$}

\address{$^{1}$ $^{2}$ Department of Mathematics, University of Illinois at Urbana--Champaign, 1409 W. Green St., Urbana, IL, 61801, U.S.A.; {aimo@illinois.edu, joe@math.uiuc.edu}}

\date{}


\abstract  
We prove that there exists an entire function for which every complex number is an asymptotic value and whose growth is arbitrarily slow subject only to the necessary  condition that the function is of infinite order.
\endabstract

\maketitle

\section{Introduction}

If $f$ is a meromorphic function in the complex plane ${\mathbb C}$, a point $w$ in the extended complex plane $ \overline{{\mathbb C}}  = {\mathbb C} \cup \{\infty\}$ is said to be an asymptotic value of $f$ if there exists a path $\gamma$ from $0$ to $\infty$ such that $f(z)$ tends to $w$ as $z$ tends to $\infty$ along $\gamma$. Iversen \cite{I1} showed that $\infty$ is an asymptotic value of every non-constant  entire function. Entire functions of finite order have at most a finite number of asymptotic values. The Denjoy--Carleman--Ahlfors Theorem states that an entire function of finite order $\rho$ has at most $2\rho$ finite asymptotic values (\cite{D1}, \cite{A}, \cite{A1}, \cite{A2}, \cite{C}). Sire \cite{S} and Iversen \cite{I2} both gave examples of entire functions of infinite order for which the set of asymptotic values has the power of the continuum. Gross \cite{G2} improved on these examples by exhibiting an entire function for which every $w$ in $ \overline{{\mathbb C}} $  is an asymptotic value.

The situation is quite different for meromorphic functions. Valiron (\cite{V1}, p.~420)  gave  examples of meromorphic functions of an arbitrary finite order, including order zero, for which the set of asymptotic values has the power of the continuum. Valiron (\cite{V2}, pp.~239--240) also showed that if $f$ is meromorphic with Nevanlinna characteristic $T(r,f)$ satisfying $T(r,f) = O( (\log r)^2 )$ as $r\to\infty$, then $f$ has at most one  asymptotic value. Eremenko \cite{E} proved that this result is sharp by showing that for every positive function $\psi(r)$, defined for $r\geq 1$, tending to $\infty$ as $r\to\infty$, there exists a meromorphic $f$ with $T(r,f) = O(\psi(r)  (\log r)^2 )$ for which each $w$ in $ \overline{{\mathbb C}} $  is an asymptotic value.

In this note we modify Gross's construction to obtain an entire function $\varphi$ having every $w\in \overline{{\mathbb C}} $ as an asymptotic value and with growth arbitrarily slow subject only to the condition that $\varphi$ has infinite order.

\begin{theorem} \label{th1}
Suppose that the function $\tilde{G} \colon [0,\infty) \to  (0,\infty)$ is continuous and strictly increasing, and that $\lim_{r\to\infty} \tilde{G} (r) = \infty$. There exist an entire function $\varphi$ and $R_0>0$ such that
\begin{equation*} 
\log\log M(r,\varphi) < \tilde{G} (r) \log r \qquad \text{for all }\,\,  r>R_0,
\end{equation*}
and such that each $w\in \overline{{\mathbb C}} $ is an asymptotic value of $\varphi$.
\end{theorem}

Here $M(r,\varphi) = \max \{ | \varphi(z) | \colon |z|=r \}$ is the usual maximum modulus function.

Our argument closely parallels that of Gross, the only substantive change being a judicious choice of the parameter $\lambda_n$ (see (12)) in order to depress the growth of $\varphi$. Other (non-essential) modifications are made in the interest of simplicity and completeness.

\section{Proof of Theorem~\ref{th1}}

Following Gross, we consider the entire function
\begin{equation*} 
\varphi_0(z) = \frac{1}{2} \left[  \frac{2}{\sqrt{\pi}} \int_0^z e^{-w^2} \, dw - e^{-z^2} + 1        \right]  .
\end{equation*}

We collect several properties of $\varphi_0$ important to our argument. We have
$$
\varphi_0(0) = 0  ,     \leqno{(1)}
$$
$$
\varphi_0(z) \hbox{ is real if } z \hbox{ is real, and}       \leqno{(2)}
$$
$$
              \qquad (i)\qquad  \lim_{x\to +\infty} \varphi_0(x) =1 ,
$$
$$
           \qquad (ii)\qquad  \lim_{x\to -\infty} \varphi_0(x) =0 ,
$$
$$
\varphi_0'(z) = e^{-z^2} \left(   z +   \frac{1}{\sqrt{\pi}}     \right)    .   \leqno{(3)}
$$
Further, 
$$
\varphi_0(x) \hbox{  is decreasing on } (-\infty, -1/\sqrt{\pi}), \hbox{ and }
$$
$$
\varphi_0(x) \hbox{  is increasing on } (-1/\sqrt{\pi} , \infty)  .
$$
$$
\hbox{ From (2) and (3) it follows that  }     \leqno{(4)}
$$
$$
              \qquad (i)\qquad  0< \varphi_0(x) <1 \qquad \hbox{ if }  x>0 ,
$$
$$
           \qquad (ii)\qquad  0> \varphi_0(x) > \frac{1}{2} (  -1-1+1  ) = -\frac{1}{2} \qquad\hbox{ if }  x<0  ,
$$
$$
M(r,\varphi_0) \leq \frac{r}{\sqrt{\pi}} e^{r^2} + \frac{1}{2}\left(  e^{r^2} +1    \right) < r e^{r^2}, \qquad r>R_0 .      \leqno{(5)}
$$
Here and later, $R_0$ refers to a large positive number, not necessarily the same at every occurrence. Furthermore, 
$$
 \qquad (i)\qquad  | \varphi_0'(z)|  \leq  e \left( 1 +  \frac{1}{\sqrt{\pi}}   \right)<4 , \qquad |z|<1  ,   \leqno{(6)}
$$
$$
 \qquad (ii)\qquad  | \varphi_0(z)|  \leq  4|z| , \qquad  |z|<1  .
$$
$$
\hbox{Suppose } 0<\varepsilon <\frac{ \pi } {4}. \,\, \hbox{From (2) and  (3) it follows that }      \leqno{(7)}
$$
$$
              \qquad (i)\qquad   \varphi_0(z)   \hbox{ converges to } 1 \hbox{ uniformly on the sector } $$
              $$ \{ r e^{i\theta} \colon r\geq 0, \, | \theta | < \frac{\pi}{4} - \varepsilon \} , $$
$$
           \qquad (ii)\qquad   \varphi_0(z)    \hbox{ converges to } 0 \hbox{ uniformly on the sector } $$
           $$\{ r e^{i\theta} \colon r\geq 0, \, | \theta - \pi | < \frac{\pi}{4} - \varepsilon \} .
$$
$$
\leqno{(8)} \qquad \text{By  (4ii) and (7ii) there exists} \,\, \alpha\in (0,\pi/2) \,\, \text{such that}
$$
$$
 |\varphi_0(z)|< \frac{3}{4} \,\, \text{in the sector}\,\, 
K := \{ re^{i\theta} \colon r\geq 0, \, |\theta - \pi |\leq \alpha/2 \} . 
$$
$$
\leqno{(9)} \qquad \text{From  (4i) and (7i) it follows that there exists a sequence} \,\, \alpha_n \,\, \text{with} \,\,
$$
$$
\alpha > \alpha_1 > \alpha_2 > \dots >0 
$$
such that
$$
\qquad (i)\qquad  | \varphi_0(z)|  \leq  2^{ \frac{1}{n} }  \quad \text{and} \quad 
| {\rm Arg}\, \varphi_0(z) | < \frac{\pi}{ 4n } 
$$
on
$$
\qquad (ii)\qquad  H_n := \{ re^{i\theta} \colon r\geq 0, \, |\theta| \leq \frac{ \alpha_n } {2} \}  .
$$

We define
$K_n := \{ re^{i\theta} \colon r\geq 0, \, |\theta -\pi | \leq \frac{ \alpha_n } {2} \} $
and note that $ | \varphi_0(z)| < \frac{3}{4}$ on $K_n\subset K$. 

We now define a family of sectors, each with  vertex at the origin, that play a central role in our argument. Let $N_n$ be a sequence of positive integers such that
$$
\leqno{(10)}  \qquad N_1=1 \,\,  \text{and}  \,\,  N_{n+1} > \frac{3\pi  } {  \alpha_n } N_n, \quad n=1,2,3,\dots  .
$$

We define two ``first generation'' sectors to be
$$
G(0) = K_1 \quad \text{and} \quad G(1)=H_1  .
$$
We note that trivially $z^{N_1}=z$ maps $G(0)$ to $K_1$ and $G(1)$ to $H_1$. 

We next define four ``second generation'' sectors as follows:

\noindent (i) $G(0,0)$ is a sector lying in $G(0)$ that $z^{N_2}$ maps onto $K_2$; 

\noindent  (ii) $G(0,1)$ is a sector lying in $G(0)$ that $z^{N_2}$ maps onto $H_2$;   

\noindent (iii) $G(1,0)$ is a sector lying in $G(1)$ that $z^{N_2}$ maps onto $K_2$;

\noindent  (iv) $G(1,1)$ is a sector lying in $G(1)$ that $z^{N_2}$ maps onto $H_2$.

Note that since each of the sectors $G(0)$ and $G(1)$ has angle opening $\alpha_1$ and since $N_2 \alpha_1 > 3\pi$ by (10), in each of the four cases above there is at least one choice of a second generation sector satisfying the stated condition. In each case, we choose one such sector and discard any others from further consideration.

Note that each second generation sector has angle opening $\theta_2$ where $N_2 \theta_2 = \alpha_2$. We further note by (10) that  $N_3 \theta_2 = N_3 \alpha_2 / N_2 > 3\pi$. Thus we may choose eight ``third generation'' sectors $G(i_1,i_2,i_3)$ with $i_j=0$ or $1$ as follows:

For all $i_1$ and $i_2$,

\noindent (i) $G(i_1,i_2,0)$ is a sector lying in $G(i_1,i_2)$ that $z^{N_3}$ maps onto $K_3$;

\noindent (ii) $G(i_1,i_2,1)$ is a sector lying in $G(i_1,i_2)$ that $z^{N_3}$ maps onto $H_3$.

As before, in each of these eight cases, we select one third generation   sector satisfying the required condition and discard all others.

We note that each third generation sector has angle opening $\theta_3$ where $N_3 \theta_3 = \alpha_3$. Since $N_4 \theta_3 = N_4\alpha_3/N_3>3\pi$ by (10), we may continue the construction to obtain $16$ fourth generation sectors, two of which lie in each third generation sector. 

We continue this construction, for each  $n$ obtaining $2^n$ sectors in the $n^{{\rm th}}$ generation, each denoted by $G(i_1,i_2,\dots ,i_n)$, $i_j=0$ or $1$, such that both $G(i_1,i_2,\dots ,i_{n-1},0)$ and $G(i_1,i_2,\dots ,i_{n-1},1)$ lie in the $(n-1)^{{\rm st}}$ generation sector $G(i_1,i_2,\dots ,i_{n-1})$. The $n^{{\rm th}}$ generation sector $G(i_1,i_2,\dots ,i_{n-1},0)$ is mapped by $z^{N_n}$  conformally onto the sector $K_n$ and $G(i_1,i_2,\dots ,i_{n-1},1)$ is mapped by $z^{N_n}$  conformally onto $H_n$. The sectors in the $n^{{\rm th}}$ generation have angle opening $\theta_n$ where $N_n \theta_n = \alpha_n$. By (10), $N_{n+1} \theta_n = N_{n+1} \alpha_n/N_n>3\pi$, allowing the construction to continue to the $(n+1)^{{\rm st}}$ generation. 

For each $n$, let $T_n$ be the union of the $2^n$ sectors of the $n^{{\rm th}}$ generation. Clearly $T_{n+1}\subset T_n$. Note that
$$
\cap_{n=1}^{\infty} T_n = \{ re^{i\theta} \colon r\geq 0, \,\, \theta\in J \}
$$
where $J\subset [0,2\pi]$ is a Cantor set.

Recall the assumptions regarding the function $\tilde{G}  \colon [0,\infty)  \to (0,\infty)$  in Theorem~\ref{th1}. Associated with the function $\tilde{G}$ there is a continuous strictly increasing function $g \colon [0,\infty)  \to (0,\infty)$ with $g(t)\to +\infty$ as $t\to +\infty$ such that $g(0)<1$ and $g(t)< \tilde{G}(t) /3$ for all $t>0$. We define an increasing unbounded sequence $r_n$ by requiring that
$$
g(r_n) = N_n 
$$
for all $n\geq 1$. By replacing $g$ with an even more slowly growing function (still denoted by $g$), if necessary, we may further specify that $r_1>e$ and $r_{n+1}>er_n$ for all $n\geq 1$, implying that
$$
\leqno{(11)}  \qquad  \log r_n >n , \quad \hbox{ for all } \, n\geq 1 .
$$

For $n\in {\mathbb N}$, let
$$
\leqno{(12)} \qquad  \lambda_n = 1/ ( 8 r_n^{g(r_n)} ) = 1/( 8 r_n^{N_n}    ) 
$$
and let
$$
\varphi_n(z) = \left(  \varphi_0(   \lambda_n z^{N_n}      )     \right)^n  .
$$

We consider the behavior of $\varphi_n$ on the $n^{{\rm th}}$ generation sector $G(i_1,i_2,\dots ,i_n)$. The function $ \lambda_n z^{N_n} $ maps $G(i_1,i_2,\dots ,i_{n-1},0)$ 
to the sector $K_n$ and consequently by (8) and (7ii)  we have
$$
\leqno{(13)} \qquad  (i)  \qquad  | \varphi_n(z) | < \left(   \frac{3}{4}   \right)^n , 
\quad z\in G(i_1,i_2,\dots ,i_{n-1},0)  ,
$$
and
$$
 \qquad (ii)  \qquad \varphi_n(z) \hbox{ tends to } 0 \hbox{ uniformly }
 \hbox{as } z\to\infty \hbox{ on } G(i_1,i_2,\dots ,i_{n-1},0)  .
$$
On the other hand, $ \lambda_n z^{N_n} $ maps the sector $G(i_1,i_2,\dots ,i_{n-1},1)$ to the sector $H_n$ and consequently by (9) and (7i) we have 
$$
\leqno{(14)} \qquad  (i)  \qquad  | \varphi_n(z) | \leq \left(  2^{ \frac{1}{n}}   \right)^n =2 , 
\quad z\in G(i_1,i_2,\dots ,i_{n-1},1) ,
$$
$$
 \qquad (ii)  \qquad | {\rm Arg}\,  \varphi_n(z) | < \frac{\pi}{4}, \quad z\in G(i_1,i_2,\dots ,i_{n-1},1)  ,
$$
and
$$
 \qquad (iii)  \qquad \varphi_n(z)  \hbox{ tends to } 1 \hbox{ uniformly }
 \hbox{as } z\to\infty \hbox{ on } G(i_1,i_2,\dots ,i_{n-1},1) .
$$

Let ${\mathbb Q}$ denote the set of rational numbers, and let
$P_1 = {\mathbb Q} \cap [0,1)$, $P_2= \{-z \colon z\in P_1\}$,
$P_3= \{iz \colon z\in P_1\}$, $P_4= \{-iz \colon z\in P_1\}$, and
$P=P_1\cup P_2\cup P_3\cup P_4$. We write $P=\{ \beta_n \colon n\in {\mathbb N} \}$ and note that $| \beta_n|<1$ for all $n$. It is elementary that given $\omega \in {\mathbb C}$, there exists $A_{\omega} \subset {\mathbb N}$ such that
$$
\leqno{(15)} \qquad  \sum_{n\in A_{\omega} } \beta_n = \omega \quad \text{and} \quad
\sum_{n\in A_{\omega} } |\beta_n| < +\infty . 
$$

We define
\begin{equation*}
\varphi(z) = \sum_{n=1}^{\infty} \beta_n \varphi_n(z) . 
\end{equation*}
We first show that $\varphi$ is entire and satisfies the required growth condition.

Suppose that $n\geq 2$ and $r_{n-1}\leq r < r_n$. If $z=re^{i\theta}$, then
$$
| \lambda_n z^{N_n} | = \lambda_n r^{N_n} = \frac{1}{8} \left(   \frac{r}{r_n}    \right)^{N_n} <  \frac{1}{8} .
$$
Thus
$$
| \varphi_0 (\lambda_n z^{N_n} ) | < \frac{1}{2} 
$$
by (6) and
$$
\leqno{(16)}  \qquad   | \beta_n \varphi_n(z) | \leq | \varphi_n(z) | < \left(   \frac{1}{2}    \right)^{n} \leq \frac{1}{4}  . 
$$

Consider $p>n$. We have
$$
| \lambda_p z^{N_p} | = \frac{1}{8} \,  \frac{  r^{N_p}   } {  r_p^{N_p}     } \leq \frac{1}{8} \left(   \frac{r_n}{r_p}    \right)^{N_p} \leq \frac{1}{8} \left(   \frac{1}{e}    \right)^{N_p}
$$
implying by (6) that
$$
| \beta_p \varphi_p(z) | \leq \left(  \frac{1}{2}  \left(   \frac{1}{e}    \right)^{N_p}     \right)^p 
$$
and consequently
$$
\leqno{(17)} \qquad   \sum_{p=n+1}^{\infty} | \beta_p \varphi_p(z) | \leq 1  .
$$
This establishes that $\varphi$ is entire.

We now consider $q$ with $1\leq q<n$. We have
$$
| \lambda_q z^{N_q} |   \leq r^{N_q} = r^{  g(r_q) } \leq  r^{  g(r) }  .
$$
By (5) we have
$$
| \varphi_0 ( \lambda_q z^{N_q} ) | \leq  r^{  g(r) } e^{ r^{  2 g(r) } }  , \quad r>R_0 ,
$$
and thus by (11)
$$
| \beta_q \varphi_q(z) | \leq r^{  q g(r) } e^{ q r^{  2 g(r) } } \leq 
r^{  (\log r) g(r)     } e^{ (\log r) r^{  2 g(r) } }  , \quad r>R_0 .
$$
Thus, again by (11), 
$$
\leqno{(18)} \qquad   \sum_{q=1}^{n-1} | \beta_q \varphi_q(z) | \leq 
(\log r)  r^{  (\log r) g(r)     } e^{ (\log r) r^{  2 g(r) } }   , \quad r>R_0 .
$$
Combining (16), (17), and (18), we have
$$
M(r,\varphi) \leq  3 (\log r) r^{   (\log r) g(r)  + r^{  2 g(r) }     }   , \quad r>R_0 .
$$
Then
$$
\log M(r,\varphi) \leq  \log 3 + \log \log r + (\log r) (   (\log r) g(r)  + r^{  2 g(r) }     )  , \quad r>R_0 .
$$
Using the inequality
$$
\log (  x+y     ) \leq \log x + \log  y  + \log 2 , \quad x\geq 1, \,\, y\geq 1, 
$$
three times, we conclude that 
\begin{eqnarray*}
\log \log M(r,\varphi) & \leq &  \log\log 3 + \log\log\log r + 2 \log\log r + \log g(r) 
\\
&{}&
+ 2 g(r) \log r + 3\log 2   , \quad r>R_0 .
\end{eqnarray*}
Then
$$
\log \log M(r,\varphi) < \tilde{G} (r) \log r  , \quad r>R_0 ,
$$
as required. 

\begin{lemma} \label{l111}
Suppose that $A\subset {\mathbb N}$ is such that
$$
\leqno{(19)} \qquad   \sum_{n\in A} | \beta_n | < +\infty 
$$
and
$$
\leqno{(20)} \qquad  \sum_{n\in A} \beta_n  = X_0 + i Y_0  .
$$
Define a sequence $i_j$ by
$$
i_j  = \left\{ \begin{matrix}  1, \quad  j\in A, \\  0,  \quad  j\notin A . \end{matrix}     \right.
$$
For   $n\in {\mathbb N}$, let $S_n$ be the associated $n^{  {\rm th} }$ generation sector
$G(i_1,i_2, \dots ,i_n)$. Let $\cap_{n=1}^{\infty} S_n = \{  re^{i\theta} \colon r\geq 0 \}$. Then
$$
\lim_{r\to \infty} \varphi( re^{i\theta}  ) = X_0+iY_0  .
$$
\end{lemma}

{\bf Proof of Lemma~\ref{l111}.} 
First observe that $S_{n+1}\subset S_n$ so that
$\cap_{n=1}^{\infty} S_n$ is a ray $L=  \{  re^{i\theta} \colon r\geq 0 \}$ for some $\theta$. Suppose that $\varepsilon>0$. There is $n_0\in {\mathbb N}$ such that
$$
\leqno{(21)} \qquad  \sum_{n\in A\atop n>n_0} |\beta_n| <   \frac{\varepsilon}{10}  ,
$$
$$
\leqno{(22)} \qquad  \left|   \sum_{n\in A\atop n\leq n_0} \beta_n  - (X_0+iY_0) \right|  <   \frac{\varepsilon}{5}  ,
$$
and
$$
\leqno{(23)} \qquad  \sum_{  n>n_0}  \left(  \frac{  3  } {   4   }       \right)^n  <   \frac{\varepsilon}{5}  .
$$

Suppose that $n\notin A$ or, equivalently, that $S_n$ has the form 
$G(i_1,i_2, \dots ,i_{n-1},0)$. Since the ray $L$ lies in $S_n$, we have from  (13i) that 
$$
\leqno{(24)} \qquad  | \varphi_n (re^{i\theta} ) | <  \left(  \frac{  3  } {   4   }       \right)^n ,
\quad r>0 ,
$$
and from (13ii) that
$$
\leqno{(25)} \qquad  \lim_{r\to\infty} \varphi_n (re^{i\theta} ) =0 .
$$
From (25) we conclude that there exists $R_1=R_1(\varepsilon, n_0)$ such that
$$
\leqno{(26)} \qquad  \sum_{n\notin A\atop n\leq n_0} |  \varphi_n (re^{i\theta} )  | <   \frac{\varepsilon}{5} , \quad r>R_1  .  
$$

Next suppose that $n\in A$ so that
$S_n=G(i_1,i_2, \dots ,i_{n-1},1)$. From (14i) and (14iii) we conclude that
$$
\leqno{(27)} \qquad   |  \varphi_n (re^{i\theta} )  | \leq 2 , \quad r>0  ,  
$$
and
$$
\leqno{(28)} \qquad  \lim_{r\to\infty} \varphi_n (re^{i\theta} ) =1 .  
$$
From (28) we conclude that there exists $R_2=R_2(\varepsilon, n_0)$ such that 
$$
\leqno{(29)} \qquad  \sum_{n\in A\atop n\leq n_0} |  \varphi_n ( re^{i\theta} )  - 1  | <   \frac{\varepsilon}{5} , \quad r>R_2  .  
$$

Now suppose that $r>\max\{ R_1,R_2\}$.  Then
\begin{eqnarray*}
&{}&
\left|   \varphi ( re^{i\theta} )  - (X_0+iY_0)      \right| 
= \left|   \sum_{n\in A\atop n\leq n_0} (  \beta_n   \varphi_n ( re^{i\theta} )      - \beta_n    ) 
\right.
\\ &{}& 
+ \left[    \sum_{n\in A\atop n\leq n_0} \beta_n     - (X_0+iY_0)    \right] 
+ \sum_{n\in A\atop n> n_0} \beta_n   \varphi_n ( re^{i\theta} ) 
\\ &{}& \left.
+ \sum_{n\notin A\atop n \leq n_0} \beta_n   \varphi_n ( re^{i\theta} )
+ \sum_{n\notin A\atop n> n_0} \beta_n   \varphi_n ( re^{i\theta} ) 
\right|  
\\ &{}&
\leq \sum_{n\in A\atop n \leq n_0} |\beta_n |  |  \varphi_n (re^{i\theta} ) - 1  | 
+ \left|  \sum_{n\in A\atop n \leq n_0}  \beta_n     - (X_0+iY_0)    \right| 
\\ &{}&
+   \sum_{n\in A\atop n > n_0} |\beta_n |  |  \varphi_n (re^{i\theta} )   | 
+ \sum_{n\notin A\atop n \leq n_0} |\beta_n |  |  \varphi_n (re^{i\theta} )   | 
\\ &{}&
+   \sum_{n\notin A\atop n > n_0} |\beta_n |  |  \varphi_n (re^{i\theta} )   | 
\\ &{}&
\leq \sum_{n\in A\atop n \leq n_0}  |  \varphi_n (re^{i\theta} ) - 1  | 
+ \left|  \sum_{n\in A\atop n \leq n_0}  \beta_n     - (X_0+iY_0)    \right| 
\\ &{}&
+ 2 \sum_{n\in A\atop n > n_0} |\beta_n | 
+ \sum_{n\notin A\atop n \leq n_0}  |  \varphi_n (re^{i\theta} )   | 
+ \sum_{n\notin A\atop n > n_0}  \left(  \frac{  3  } {   4   }       \right)^n
\\ &{}&
< \frac{\varepsilon}{5} + \frac{\varepsilon}{5} + \frac{\varepsilon}{5} + \frac{\varepsilon}{5} + \frac{\varepsilon}{5} = \varepsilon ,
\end{eqnarray*}
where we have used (27), (24), (29), (22), (21), (26), and (23).

This proves Lemma~\ref{l111}. 

Lemma~\ref{l111} implies that the conclusion of Theorem~\ref{th1} holds for $\omega = X_0+iY_0\in {\mathbb C}$, as (15) tells us that there is a set $A \subset {\mathbb N}$ satisfying (19) and (20). 

Recall that by Iversen's Theorem \cite{I1} from 1914, infinity is an asymptotic value for every non-constant entire function, hence also for $\varphi$. We further remark that in this particular case,  there is an asymptotic path for the asymptotic value $\infty$ that is a ray. Let $A \subset {\mathbb N}$ be such that $\beta_n>0$ for all $n\in A$ and $\sum_{n\in A} \beta_n=+\infty$. Let $i_j=1$ if $j\in A$ and $i_j=0$ if $j\notin A$. Let $S_n = G(i_1,i_2, \dots ,i_n)$ be the associated $n^{  {\rm th} }$ generation sector. Note that $S_{n+1}\subset S_n$. Let $L=\cap_{n=1}^{\infty} S_n$ be the ray $L=\{ re^{i\theta} \colon r\geq 0 \}$. 

Suppose that $M>3$. Let $n_0 \in A$ be such that 
$\sum_{n\in A\atop n \leq n_0} \beta_n > 4 M $. 
Note that
$$
 \lim_{r\to\infty} \varphi_n (re^{i\theta} ) = 1  
$$
for all $n\in A$ by (14iii) since $L\subset S_n$. Thus there exists $R_0>1$ such that
$$
{\rm Re}\, \sum_{n\in A\atop n \leq n_0}  \beta_n   \varphi_n ( re^{i\theta} ) 
> \frac{1}{2} \sum_{n\in A\atop n \leq n_0}  \beta_n  > 2 M , \quad r>R_0 .
$$
By (14ii), ${\rm Re}\,  \beta_n   \varphi_n ( z ) \geq 0$ for all $n\in A$ and $z\in S_n$. In particular, ${\rm Re}\,  \beta_n   \varphi_n ( re^{i\theta} ) \geq 0$ for all $n\in A$ and all $re^{i\theta}\in L$. Thus for all $r>R_0$ we have using (24) 
\begin{eqnarray*}
&{}&
{\rm Re}\, \varphi(    re^{i\theta}     ) 
= {\rm Re}\, \sum_{n\in A\atop n \leq n_0}   \beta_n   \varphi_n ( re^{i\theta} )
+ {\rm Re}\, \sum_{n\in A\atop n > n_0}   \beta_n   \varphi_n ( re^{i\theta} )
+ {\rm Re}\, \sum_{n\notin A}   \beta_n   \varphi_n ( re^{i\theta} )
\\ &{}&
> 2 M + 0 - \sum_{n=1}^{\infty} \left(  \frac{  3  } {   4   }       \right)^n
= 2 M - 3 > M  .
\end{eqnarray*}
Thus $\varphi(z)$ tends to $\infty$ on $L$. This completes the proof of Theorem~\ref{th1}. 
 
\section{Concluding remarks}

1. For every complex number $X_0+iY_0$, the set of rays $\{ re^{i\theta} \colon r\geq 0, \, \theta \in J \}$ on which $\varphi$ tends to $X_0+iY_0$ is large in the sense that $J$ has the power of the continuum. Let $B$ be any infinite subset of ${\mathbb N}$ such that $\beta_n>0$ for all $n\in B$ and $\sum_{n\in B} \beta_n=1$. Let $C$ be an arbitrary subset of $B$. By an obvious modification of (15), there exists a set $A_C \subset {\mathbb N} \setminus B$ such that 
$$
\sum_{ n\in A_C } |\beta_n | < +\infty ,\quad \text{and} \quad \sum_{ n\in A_C } \beta_n  =X_0+iY_0
-  \sum_{ n\in C } \beta_n  . 
$$
Set $A=C \cup A_C$. Then $A$ satisfies conditions (19) and (20) of Lemma~\ref{l111}. Also note that $A\cap B=C$. Thus there is an injection from the collection of all subsets $C$ of $B$ into the collection of all subsets $A$ of ${\mathbb N}$ satisfying (19) and (20). Lemma~\ref{l111} establishes an injection from the collection of all such subsets $A$ of ${\mathbb N}$ into $J$. Since the collection of all subsets $C$ of $B$ has the power of the continuum, we conclude that the set $J$ does as well. A similar argument applies to the set of rays on which $\varphi$ tends to infinity. 

2. Suppose that $m \in {\mathbb N}$ and let 
$\tilde{S}_m = G(i_1,i_2,\dots ,i_m)$ be any sector of the $m^{  {\rm th} }$ generation. 
Suppose that $X_0+iY_0\in {\mathbb C}$. We claim that there exists a ray $L\subset \tilde{S}_m$ with $ L = \{ re^{i\theta} \colon r\geq 0  \}$ such that 
$\lim_{r\to\infty} \varphi (re^{i\theta} ) =X_0+iY_0$. 

To justify this claim, let
$A^* = \{ j \colon 1\leq j\leq m, \, i_j=1 \}$. 
Using a minor modification of (15), let
$A^{**} \subset {\mathbb N} \setminus \{1,2,\dots ,m\}$ be such that
$$
\sum_{n\in A^{**} } |\beta_n| <+\infty \quad \text{and}\quad 
\sum_{n\in A^{**} } \beta_n = X_0+iY_0 - \sum_{n\in A^{*} } \beta_n .
$$
Set $A= A^* \cup A^{**}$. Note that $A$ satisfies (19) and (20). For $j\in {\mathbb N}$, let $i_j=1$ if $j\in A$ and $i_j=0$ if $j\notin A$. Let $S_n = G(i_1,i_2, \dots ,i_n)$ be the associated $n^{  {\rm th} }$ generation sector.
Note that $S_{n+1}\subset S_n$ and $S_m= \tilde{S}_m$. 
Let $L=\cap_{n=1}^{\infty} S_n = \{ re^{i\theta} \colon r\geq 0  \}$. 
Lemma~\ref{l111} asserts that 
$\lim_{r\to\infty} \varphi (re^{i\theta} ) =X_0+iY_0$. 
Since $L\subset S_m= \tilde{S}_m$, our claim is established.
It can similarly be shown that $\varphi$ tends to infinity on some ray in $ \tilde{S}_m$.

\end{document}